\documentstyle[12pt,amsfonts,amscd]{amsart}
\newtheorem{Theorem}{Theorem}[section]
\newtheorem{Definition}[Theorem]{Definition}
\newtheorem{Proposition}[Theorem]{Proposition}

\newtheorem{Lemma}[Theorem]{Lemma}
\newtheorem{Corollary}[Theorem]{Corollary}
\theoremstyle{remark}
\newtheorem{Example}[Theorem]{Example}

\def\CC{{\Bbb C}}
\def\RR{{\Bbb R}}
\def\ZZ{{\Bbb Z}}

\def\be{\begin{enumerate}}
\def\ee{\end{enumerate}}
\def\bT{\begin{Theorem}}
\def\eT{\end{Theorem}}
\def\bP{\begin{Proposition}}
\def\eP{\end{Proposition}}
\def\bD{\begin{Definition}}
\def\eD{\end{Definition}}
\def\bE{\begin{Example}}
\def\eE{\end{Example}}
\def\bL{\begin{Lemma}}
\def\eL{\end{Lemma}}
\def\bC{\begin{Corollary}}
\def\eC{\end{Corollary}}

\begin{document}
\title{Osgood-Hartogs type properties of power series and smooth functions}
\author{Buma L. Fridman and  Daowei Ma}
\begin{abstract} We study the convergence of a formal power series of two variables if
its restrictions on curves belonging to a certain family are convergent. Also
analyticity of a given $C^\infty $ function $f$ is proved when the
restriction of $f$ on analytic curves belonging to some family is analytic. Our results generalize two known statements: a theorem of P. Lelong and the Bochnak-Siciak Theorem. The questions we study fall into the category of ``Osgood-Hartogs-type'' problems. \end{abstract}


\address{ buma.fridman@@wichita.edu, Department of Mathematics,
Wichita State University, Wichita, KS 67260-0033, USA}
\address{ dma@@math.wichita.edu, Department of Mathematics,
Wichita State University, Wichita, KS 67260-0033, USA}
\maketitle \setcounter{section}{-1}
\section{Introduction}
The following Hartogs theorem is a fundamental result in complex analysis: a
function $f$ in ${\Bbb C}^n$, $n>1$, is holomorphic if it is holomorphic in each
variable separately, that is, $f$ is holomorphic in ${\Bbb C}^n$ if for each axis it is
holomorphic on every complex line parallel to this axis. In the last
interpretation this statement leads to a number of questions summarized by K. Spallek, P. Tworzewski, T. Winiarski
\cite{ST} the following way: ``Osgood-Hartogs-type problems ask for properties of
`objects' whose restrictions to certain `test-sets' are well known''.
\cite{ST} has a number of examples of such problems. Here are two classical examples.

P.~Lelong's theorem \cite{Le}. A formal power series $g(x,y)$ converges in some neighborhood of the origin if there exists a set $E\subset \CC$ of positive capacity, such that for each $s\in E$, the formal power series $g(x,sx)$ converges in some neighborhood of the origin (of a size possibly depending on $s$).

The Bochnak-Siciak Theorem [Bo,Si]. Let $f\in C^\infty (D),$ $D$ is a
domain, $0\in D\subset {\Bbb R}^n.$ Suppose  $f$ is analytic on every line segment
through $0$. Then $f$ is analytic in a neighborhood of $0$ (as a function
of $n$ variables).

In many  articles the same two ``objects'' are usually considered: power series and functions of
several variables. The ``test-sets'' in many cases form a family of linear
subspaces of lower dimension. For example, articles by S.~S.~Abhyankar, T.~T.~Moh \cite{AM}, N.~Levenberg and R.~E.~Molzon, \cite{LM}, R.~Ree \cite{Re}, A.~Sathaye \cite{Sa}, M.~A.~Zorn \cite{Zo} and others consider the convergence of formal power series of several
variables provided the restriction of such a series on each element of a
sufficiently large family of linear subspaces is convergent. T.S.~Neelon \cite{Ne2,Ne3} proves that a formal power series is convergent if its restrictions to certain families of curves or surfaces parametrized by polynomial maps are convergent. Articles by J.~Bochnak \cite{Bo}, T.~S.~Neelon \cite{Ne,Ne2}, J.~Siciak \cite{Si}, and others
prove that a function of several variables is highly smooth (or even
analytic) if it is smooth enough on each of a sufficiently large set of linear or algebraic curves (or surfaces of lower dimension). The publication by E.~Bierstone, P.~D.~Milman, A.~Parusinski \cite{BM} provides an interesting example of a non-continuous function in ${\Bbb R}^2$
that is analytic on every analytic curve.

In this article we also consider both: power series with complex coefficients and
functions in a neighborhood of the origin in ${\Bbb R}^2$. As far as
``test-sets'' we consider separately two families. They are derived the
following way. First consider a non-linear analytic curve $\Gamma =\{x,\gamma (x)\},$ $%
\gamma (0)=0$. One family $\Im _1$ is a set of dilations of $\Gamma$: $\Im
_1=\{sx,s\gamma (x)\},s\in \Lambda _1\}$, where $\Lambda _1\subset {\Bbb R}$
is a closed subset of $\CC$ of positive capacity; the other family $\Im _2$ consists
of curves $\Gamma _{\theta \text{ }}$ $(\theta \in \Lambda _2)$ each of which is a
rotation of $\Gamma $ about the origin by an angle $\theta $; $\Lambda _2$
is a subset of $[0,2\pi ]$ of positive capacity. If $f\in C^\infty $
and its restriction on every curve of $\Im _1$ ($\Im _2$ respectively) can
be extended as an analytic function in a neighborhood of that curve, then $f$ 
is real analytic in a neighborhood of the region covered by the curves of $\Im _1$ ($\Im _2$
respectively). (For precise statements see Theorems 2.1, 2.2).

We start however with two results related to power series. First we prove a generalization of P.~Lelong's theorem.
Namely, if $g(x,y)$ is a formal power series and $h(x)$, $h(0)=0$, is a convergent
power series such that the inhomogeneous dilations $g(s^\sigma x,s^\tau h(x))
$ are convergent for sufficiently many $s$ ($\sigma ,\tau $ are fixed), then 
$g(x,y)$ is convergent (for the precise statement see
Theorem~1.1). Theorem 1.2 is devoted to a reverse claim: if $h(x)$ is
a formal power series and $g(s^\sigma x,s^\tau h(x))$ converges for sufficiently many $s$,
then $h(x)$ is convergent (see Theorem 1.2 for exact statement).

The results in this paper do not carry over in a routine way to dimensions greater than two. We intend to study corresponding problems for higher dimensions in future work.

\section{On the convergence of a power series in two variables}

Let $\CC[[x_1, x_2, \dots, x_n]]$ denote the set of (formal) power series 
$$g(x_1,\dots, x_n)=\sum_{k_1,\dots, k_n\ge0} a_{k_1\dots k_n} x_1^{k_1}\cdots x_n^{k_n} $$
of $n$ variables with complex coefficients. Let $g(0)=g(0,\dots,0)$ denote the coefficient $a_{0,\dots,0}$. A  power series equals 0 if all of its coefficients $a_{k_1\dots k_n}$ are equal to 0. A  power series $g\in \CC[[x_1, x_2, \dots, x_n]]$ is said to be convergent if there is a constant $C=C_g$ such that $|a_{k_1\dots k_n} |\le C^{k_1+\cdots +k_n}$ for all $(k_1,\dots,k_n)\not=(0,\dots,0)$. If $g$ is convergent, then it represents a holomorphic function in some neighborhood of $0$ in $\CC^n$. If $g\in \CC[[x_1, x_2, \dots, x_n]]$ and $s\in \CC^n$, then $g_s(t):= g(s_1t,\dots, s_nt)$ is well defined and belongs to $\CC[[t]]$. By \cite{Zo}, $g$ is convergent  if and only if $g_s(t)$ is convergent for each $s\in\CC^n$.  The partial derivatives of a power series are well defined even when it is divergent (not convergent). For example, if $g\in \CC[[x,y]]$ and if $g=\sum a_{ij}x^iy^j$, then
$$g'_y=\frac{\partial g}{\partial y}=\sum ja_{ij}x^iy^{j-1}.$$ Thus $g'_y\not =0$ simply means that $g\not\in \CC[[x]]$.
If $g\in\CC [[x,y]]$,  and if $h\in\CC [[x]]$ with $h(0)=0$,  then $g(x, h(x))$ is a well-defined element of $\CC[[x]]$. 

As mentioned above, a  lot of work has been done on the convergence of a power series with the assumption that the series is convergent after restriction to sufficiently many subspaces (see \cite{AM,LM,Le,Si,Si2}). 

We consider substitution of a power series $y=h(x)$ into an inhomogeneous dilation $g(s^\sigma x, s^\tau y)$ of a series $g(x,y)$, where $\sigma, \tau$ are integers. 

Let $$Q:=\{(\sigma,\tau): \sigma, \tau\in \ZZ, (\sigma,\tau)\not=(0,0)\}.$$  Let $cap(E)$ denote the (logarithmic) capacity of a closed set $E$ in the complex plane $\CC$.

We now present our two main theorems. 

\bT \label{MT1} Let $g\in\CC [[x,y]]$ be a power series of two variables $x, y$, let $h\in \CC[[x]]$ be a non-zero convergent power series with $h(0)=0$, let $E$ be a closed subset of ${\Bbb{C}}\backslash \{0\}$ with $cap(E)>0$, and let $(\sigma,\tau) $ be a pair in the set $Q$. Assume, in case $\sigma\tau>0$, that $h(x)$ is not a monomial of the form $b_kx^k$ with $\sigma k-\tau=0$. Suppose that $g(s^\sigma x, s^\tau h(x))$ is convergent for each $s\in E$. Then $g$ is convergent.
 \eT

\bT \label{MT2} Let $g\in\CC [[x,y]]$ be a power series with $g'_y\not=0$, let $h\in \CC[[x]]$ be a non-zero power series with $h(0)=0$, let $E$ be a closed subset of ${\Bbb{C}}\backslash \{0\}$ with $cap(E)>0$, and let $(\sigma,\tau) $ be a pair in the set $Q$ with $\sigma\tau>0$. Suppose that $g(s^\sigma x, s^\tau h(x))$ is convergent for each $s\in E$. Then $h$ is convergent.
 \eT

The examples in Section~\ref{exam} show that in the above two theorems, if any condition is dispensed with, then the resulting statement is false. We now prove some auxiliary results.

The following theorem is a consequence of a result by B. Malgrange \cite{Ma}.
We present below an independent short proof of the statement.

\bT \label{L1} Let $g\in \CC[[x_1,\dots,x_n,y]]$ with $g'_y\not =0$, and let $h  \in \CC[[x_1,\dots,x_n]]$ with $h(0 )=0$. Suppose that $g$ and $g(x_1,\dots,x_n, h(x_1,\dots, x_n))$ are convergent. Then $h$ must be convergent.
\eT
\begin{pf}Let $f\in \CC[[x_1,\dots,x_n,y]]$ be defined by 
$$f(x_1,\dots,x_n,y)=g(x_1,\dots,x_n,y)-g(x_1,\dots,x_n,h(x_1,\dots,x_n)).$$ Then $f$ is convergent and  $f(x_1,\dots,x_n,h(x_1,\dots,x_n))=0$. Fix an $s=(s_1,\dots,s_n)\in\CC^n$. Let $f_s(t,y)\in\CC[[t,y]]$ be defined by $f_s(t,y)=f(s_1t,\dots,s_nt,y)$. Then $f_s(t,y)$ is convergent and $f_s(t,h_s(t))=0$. By Weierstrass preparation theorem (see, {\it e.g.}, \cite{GH}, p.~8), there is a nonnegative integer $k$ such that \ $f_ s(t,y)={t^k}P(t,y)Q(t,y)$, where $P(t,y)=y^m+a_1(t)y^{m-1}+\cdots+a_m(t)$ is a polynomial in $y$ with coefficients being convergent power series in $t$, and $Q(t,y)$ is a convergent power series with $Q(0,0)\not=0$. Hence $P(t, h_s(t))=0$. It is known (see \cite{Fu},Theorem 4.12, p.~73) that there is a positive integer $r$ such that $P(t^r,y)$ splits into linear factors in $y$: $$P(t^r, y)=(y-u_1(t))\cdots (y-u_m(t)),$$ where the $u_j(t)$ are convergent power series. Thus $$0=P(t^r,h_s(t^r))=(h_s(t^r)-u_1(t))\cdots(h_s(t^r)-u_m(t)).$$ It follows that $h_s(t^r)=u_j(t)$ for some $j$. Therefore $h_s(t)$ is convergent. Since $h_s(t)$ is convergent for each $s\in\CC^n$, the series $h(x_1,\dots,x_n)$ must be  convergent.
\end{pf}

Let $E$ be a closed bounded set in the complex plane. The transfinite diameter of $E$ is defined as
$$d_{\infty}(E)=\lim_{n}(\max\{\Pi_{ i<j}|z_i-z_j|^{2/n(n-1)}:z_1,\dots,z_n\in E\}) .$$  For a probability measure $\mu$ on the compact set $E$, the logarithmic potential of $\mu$ is $$p_\mu(z)=\lim_{N\to\infty}\int\min(N,\log{1\over|z-\zeta|})\,d\mu(\zeta),$$ and the capacity of $E$ is defined by $$cap(E)=\exp(-\min_{\mu(E)=1}\sup_{z\in\CC} p_{\mu}(z)).$$ It turns out that $d_{\infty}(E)=cap(E)$ (see \cite{Ah2}, pp.~23--28). It follows from the definition of the transfinite diameter that $$cap(E)=\lim(cap(E_n)) \;\mbox{if}\;E=\cap E_n,\;E_1\supset E_2\supset \cdots,$$ and from the definition of the capacity that \begin{equation}\label{union}cap(E)=\lim(cap(E_n)) \;\;\mbox{if}\;E=\cup E_n,\;\;E_1\subset E_2\subset \cdots.\end{equation}

If $E$ is a closed set, its capacity can be defined by 
$$cap(E)=\lim_n cap(E\cap \{|x|\le n\}).$$

\bL \label{L3} {\rm (Bernstein Inequality.)} Let $E$ be a compact set in the complex plane with $cap(E)>0$. 
Then there exists a positive constant $C=C_E$, depending only on $E$,  such that for each positive integer $n$ and each polynomial $P(z)=\sum a_k z^k \in\CC[z]$ of degree $n$, each coefficient $a_k$, $0\le k\le n$, of $P(z)$ satisfies
$$|a_k|\le C^n \max_{z\in E}|P(z)|.$$ \eL
Proposition 4.6 in \cite{Ne2} can be used to prove this statement. Also (we thank Nessim Sibony 
for pointing this out to us) this Lemma follows from considerations in 
N. Sibony, Sur la fronti\`{e}re de Shilov des domaines de ${\Bbb C}^n$, Math. Annalen 273 (1985), 115-121. We present here an independent short proof of this Lemma.
\begin{pf} Without loss of generality we assume that $\max_{z\in E}|P(z)|=1$. Let $\Omega$ be the unbounded component of the complement of $E$ in $\CC$. It is known that $\Omega$ has a Green's function with a pole at $\infty$ (see \cite{Ah}, and \cite{Ah2}, pp.~25--27). The Green's function is harmonic in $\Omega$, 0 on $\partial \Omega$, and  its asymptotic behavior at $\infty$ is of the form $$u(z)=\log|z|-\log\alpha+o(1),$$
where $\alpha:=cap(E)$. On applying the maximum principle to the subharmonic function
\newline$\log|P(z)|-(n+\epsilon)u(z)$, we obtain $|P(z)|\le e^{nu(z)}$ for $z\in \Omega$.
Choose an $R>1$ so that $E\subset \{z: |z|<R\}$. Set $C=\max_{|z|=R}e^{u(z)}$. Then $|P(z)|\le C^n$ if $|z|=R$, and $$|a_k|=|{1\over 2\pi i}\int_{|z|=R}{P(z)\over z^{k+1}}\,dz|\le R^{-k}\max_{|z|=R}|P(z)|\le C^n.$$ This proves the lemma. \end{pf}

%\bT \label{PS} Let $g=\sum a_{ij}x^iy^j \in\CC [[x,y]]$ be a formal power series in $x,y$ with $g'_y(x,y)\not=0$, and let $h=\sum_k b_k x^k\in \CC[[x]]$ be a nonzero power series in $x$ with $b_0=0$. Suppose that $h$ is nonlinear (i.e., $h\not=b_1x$) and that $f(x, h_s(x))$ converges for each $s$ in some interval $[\alpha, \beta]$ of positive length. Then $f$ converges and $h$ converges.  \eT

{\it Proof of Theorem~\ref{MT1}. } We assume that $a_{00}=g(0,0)=0$, $E$ is bounded, $gcd(\sigma,\tau)=1$, $\sigma\ge 0$,  and, in case $\sigma=0$, $\tau=-1$.  This does not cause any loss of generality. Indeed, if $E$ is unbounded, we set $E_n=\{s\in E: n\ge |s|\ge 1/n\}$. Since $\lim cap(E_n)=cap(E)>0$, the set $E_n$ has positive capacity when $n$ is sufficiently large. On replacing $E$ by $E_n$, we obtain that $0\not\in E$ and $E$ is bounded. If $d:=gcd(\sigma,\tau)>1$, we can replace $(\sigma,\tau)$ by $(\sigma/d,\tau/d)$, and $E$ by the set $\{s\in\CC: s^d\in E\}$. Finally, if $\sigma <0$, or if $(\sigma, \tau)=(0,1)$, we can replace $(\sigma,\tau)$ by $(-\sigma,-\tau)$, and $E$ by $\{s\in\CC: s^{-1}\in E\}$.

 Let $$h(x)=\sum_{i=1}^\infty b_{i}x^i.$$ Then 
$$h(x)^j=\sum_{k=j}^\infty c_{jk}x^k,$$ where
$$c_{jk}=\sum_{l_1+\cdots+l_j=k}b_{l_1}\cdots  b_{l_j}.$$  Note that $c_{jk}=0$ for $k<j$. Hence
$$g(s^\sigma x,s^\tau h(x))=\sum_{i,j,k}a_{ij}c_{jk}s^{\sigma i+\tau j}x^{i+k }=\sum_{p=1}^{\infty}(\sum_{q=-\tau ^-p}^{(\sigma+\tau^+)p} d_{pq}s^q)x^p,$$ where $\tau^+=\max(0,\tau)$, $\tau^-=-\min(0,\tau)$, and
\begin{equation}\label{dpq}d_{pq}=\sum_{\sigma i+\tau j=q}a_{ij}c_{j,p-i}.\end{equation}

Note that for each $p\ge 1$ and each $q\in\ZZ$, the sum (\ref{dpq}) contains only a finite number of non-zero terms. Let $u_p(s)=\sum_qd_{pq}s^q$. Then $s^{\tau^-p}u_p(s)$ is a polynomial in $s$ of degree at most $(\sigma+|\tau|)p$,  and $g(s^\sigma x,s^\tau h(x))=\sum u_p(s)x^p$. For $s\in E$, since $g(s^\sigma x,s^\tau h(x))$ is convergent, its coefficients $u_p(s)$ satisfy  $|u_p(s)|\le C_s^p$ for some positive constant $C_s$, possibly depending on $s$, and $p=1, 2,\dots$.  Set, for $n=1,2,\dots$, 
$$E_n=\{s\in E: |u_p(s)|\le n^p\;\; \forall p>0\}.$$ The sequence $(E_n)$ is an increasing sequence of closed sets. Since $\lim cap(E_n)=cap(E)>0$, the set $E_n$ has positive capacity for some $n$. On replacing $E$ by $E_n$, we obtain $|u_p(s)|\le n^p$ for $s\in E$ and $p=1, 2, \dots$. The polynomial $s^{\tau^-p}u_p(s)$ is of  degree at most $(\sigma+|\tau|)p$, and satisfies $$ |s^{\tau^-p}u_p(s)|\le M^{\tau^-p}n^p,\;\; s\in E, $$ where $M=\max_E |s|$.  By Lemma~\ref{L3}, the coefficients of the above mentioned polynomial satisfy $|d_{pq}|\le C_E^{(\sigma+|\tau|)p}M^{\tau^-p}n^p$, where $C_E$ is the constant in Lemma~\ref{L3}, depending only on $E$. Set $C=C_E^{\sigma+|\tau|}M^{\tau^-}n$. Then \begin{equation}\label{est}|d_{pq}|\le C^p.\end{equation}
 
Let \begin{equation}\label{gq}g_q(x,y)=\sum_{\sigma i+\tau j=q}a_{ij}x^iy^j,\end{equation} and let $\phi_q(x)=g_q(x,h(x))$,  for $q\in \ZZ$. Then $g_q\in \CC[[x,y]]$ in general, and  it is a polynomial when $\sigma,\tau>0$. It is straightforward to verify that 
\begin{equation}\label{conv}\phi_q(x)=g_q(x,h(x))=\sum_{p=1}^{\infty} d_{pq}x^p.\end{equation} The series $\phi_q(x)$ is convergent because of (\ref{est}). Choose a positive number $r<1/C$, where $C$ is the constant in (\ref{est}), so that $h(x)$ converges in a neighborhood of the closed ball $\{x\in\CC: |x|\le r\}$ and $h(x)\not=0$ when $0<|x|\le r$. Let $m=\min_{|x|=r}|h(x)|$. Then $m>0$. For $x\in\CC$, $|x|\le r$,  $$|\phi_q(x)|\le \sum |d_{pq}| |x|^p\le \sum (Cr)^p={1\over 1-Cr}.$$

We now consider two cases: (i) $\sigma\tau>0$, and (ii) $\sigma\tau\le0$.

Case (i). $\sigma>0, \tau>0$. Let 
\begin{equation}\label{omq}\Omega_q=\{(i,j): i,j\in\ZZ,\, i,j\ge0,\,\sigma i+\tau j=q\}.\end{equation} Let $\omega_q$ be the cardinality of $\Omega_q$. It is clear that $\omega_q\le q+1$. Fix a $q\ge1$ so that $\omega_q>0$. Let $(\lambda, \mu)$ be the element of $\Omega_q$ so that $\mu$ is the minimum. Then $$\Omega_q={\{(\lambda-k\tau,\,\mu+k\sigma): k=0,1,\dots,\omega_q-1\}},$$ and $$g_q(x,y)=x^{\lambda}y^{\mu}\sum_{k=0}^{\omega_q-1} a_{\lambda-k\tau,\,\mu+k\sigma}(x^{-\tau} y^\sigma)^{k}.$$ Let $$\psi_q(t)=\sum_{k=0}^{\omega_q-1} a_{\lambda-k\tau,\,\mu+k\sigma}t^{k},$$ so that $g_q(x,y)=x^{\lambda}y^{\mu}\psi_q(x^{-\tau} y^\sigma)$, and \begin{equation}\label{psiest}\psi_q(x^{-\tau} h(x)^\sigma)=x^{-\lambda}h(x)^{-\mu}\phi_q(x).\end{equation}

Let $u(x)=x^{-\tau}h(x)^\sigma$, $S=\{x\in\CC:|x|=r\}$, and $F=u(S)$. Since $h(x)$ is not a monomial of the form $b_kx^k$ with $\sigma k-\tau=0$, the function $u(x)$ is a non-constant meromorphic function, hence $F$ has positive capacity. For $t=x^{-\tau}h(x)^\sigma \in F$, we obtain, by (\ref{psiest}), that
\begin{equation}\label{psi2}|\psi_q(t)|\le {r^{-\lambda}m^{-\mu}\over 1-Cr}\le {(1+r^{-1}+m^{-1})^{\lambda+\mu}\over 1-Cr}.\end{equation} 
The summand $1$ in the right-hand side of the above inequality is included to ensure that the numerator is greater than 1 as needed later. 
Hence $|\psi_q(t)|\le L^q$ on $F$, where 
$$L={1+r^{-1}+m^{-1}\over 1-Cr},$$ for $\lambda+\mu\le q$. By Lemma~\ref{L3}, the coefficients of $\psi_q$ are bounded by $L^qC_F^{\omega_q-1}$. Thus for $(i,j)\in\Omega_q$, $$|a_{ij}|\le L^qC_F^{\omega_q-1}\le (L+C_F)^{2q}\le (L+C_F)^{2(\sigma+\tau)(i+j)},$$ or $|a_{ij}|\le K^{i+j}$, where $K=(L+C_F)^{2(\sigma+\tau)}$. Note that the number $K$ does not depend on $q$. It follows that $$|a_{ij}|\le K^{i+j}, \;\;\mbox{if} \;\; \sigma i+\tau j\ge 1.$$ This proves that $g$ is convergent.

Case (ii). $\sigma\ge 0$, $\tau\le 0$. In this case the set $\Omega_q$ in (\ref{omq}) can be written as
$$\Omega_q={\{(\lambda-k\tau,\,\mu+k\sigma): k=0,1,2,\dots\}},$$ where $(\lambda,\mu)$ is the element in $\Omega_q$ with least value of $\mu$ when $\sigma>0$, and $(\lambda,\mu)=(0, -q)$ when $(\sigma,\tau)=(0,-1)$. Let 
$$\psi_q(t)=\sum_{k=0}^{\infty} a_{\lambda+k|\tau|,\,\mu+k\sigma}t^{k}.$$ Then $g_q(x,y)=x^\lambda y^\mu\psi_q(x^{|\tau|}y^\sigma)$. The formal power series $\psi_q(t)$ satisfies $\phi_q(x)=x^\lambda h(x)^\mu\psi_q(x^{|\tau|}h(x)^\sigma)$. Since $x^\lambda h(x)^\mu$ and $\phi_q(x)$ are convergent, $\alpha(x):=\psi_q(x^{|\tau|}h(x)^\sigma)$ has to be convergent. Write $x^{|\tau|}h(x)^\sigma= c x^\nu+\cdots$, $c\not=0$.  There is a power series $\beta(x)$, also convergent in a neighborhood of $\{|x|\le r\}$, such that $x^{|\tau|}h(x)^\sigma=\beta(x)^\nu$. Reducing $r$ if necessary, we assume that $\beta(x)$ is univalent in a neighborhood of $\{|x|\le r\}$. Note that the reduction in the value of $r$ is independent of $q$. The set $\{\beta(x): |x|<r\}$ contains an open disc $\{z\in \CC: |z|<\delta\}$. The series $\beta(x)$ has an inverse $\gamma(z)$, convergent in $\{z\in \CC: |z|<\delta\}$, such that $\gamma(\beta(x))=x$ and $\beta(\gamma(z))=z$. Now $\psi_q(z^\nu)$ is converge
 
 nt in $\{|z|<\delta\}$, so $\psi_q(t)$ is convergent in $\{|t|<\delta^\nu\}$.  Let  $t\in\CC$ with $|t|<\delta^\nu$. Then $t=z^\nu$ for some $z$ with $|z|<\delta$, and $z=\beta(x)$ for some $x$ with $|x|<r$. Hence $$|\psi_q(t)|=|\psi_q(\beta(x)^\nu)|=|\alpha(x)|\le\max_{|x|=r}|\alpha(x)|.$$ Thus $$\sup_{|t|<\delta^\nu}|\psi_q(t)|\le \max_{|x|=r}|{\phi_q(x)\over x^\lambda h(x)^\mu}|\le{r^{-\lambda}m^{-\mu}\over 1-Cr}.$$ By the Cauchy estimates, the coefficients of $\psi_q$ satisfy $$|a_{\lambda+k|\tau|,\,\mu+k\sigma}|\le {r^{-\lambda}m^{-\mu}\over 1-Cr}\delta^{-k\nu}\le {(1+r^{-1}+m^{-1}+\delta^{-\nu})^{\lambda+\mu+k}\over1-Cr}.$$ The summand $1$ in the right-hand side of the above inequality is included to ensure that the numerator is greater than 1 as needed later.
It follows that, for $(i,j)\in \Omega_q$, $$|a_{ij}|\le ({1+r^{-1}+m^{-1}+\delta^{-\nu}\over1-Cr})^{i+j}.$$ Note that the number $K:=(1+r^{-1}+m^{-1}+\delta^{-\nu})/(1-Cr)$ does not depend on $q$. Therefore, $|a_{ij}|\le K^{i+j}$ for all $(i,j)$. This proves that $g$ is convergent. \qed

\vspace{5pt}

{\it Proof of Theorem~\ref{MT2}.  } This proof and the proof of Theorem~1.1 share the discussion through equation~(5). Note that the convergence of $h$ has not been used in the derivation of (5). We define polynomials $g_q(x,y)$ by (\ref{gq}). Then $g_q(x,h(x))$ are convergent by (\ref{est}) and (\ref{conv}). 
Since $g'_y(x,y)\not=0$, $\partial g_q/\partial y\not=0$ for some $q$. It follows from Theorem~\ref{L1} that $h(x)$ is convergent.  \qed

\vspace{5pt}

For $h\in \CC[[x]]$ with $h(0)=0$, let $h_s(x)=s^{-1}h(sx)$.
\bC \label{Dil}Let $g\in\CC [[x,y]]$ be a power series, let $h\in \CC[[x]]$ be a non-zero and non-linear power series with $h(0)=0$, and let $E$ be a closed subset of ${\Bbb{R}}\backslash \{0\}$  with $cap(E)>0$. Suppose that $g(x, h_s(x))$ is convergent for each $s\in E$. Then $g$ is convergent. \eC
\begin{pf} If $g'_y=0$ then the statement holds. Suppose $g'_y\not=0$. For $s\not=0$, $g(x, h_s(x))$ is convergent if and only if $g(s^{-1} x, h_s(s^{-1}x))=g(s^{-1}x,s^{-1}h(x))$ is convergent. By Theorem~\ref{MT2}, $h$ is convergent. Then $g$ is convergent by Theorem~\ref{MT1}.\end{pf}

For $f\in \CC[[x,y]]$ and $\theta\in[0,2\pi]$, write $$f_\theta(x,y)=f(x\cos\theta-y\sin\theta, x\sin\theta+y\cos\theta).$$

\bT \label{Rot}Let $f\in\CC [[x,y]]$ be a power series, let $h\in \CC[[x]]$ be a convergent power series with $h(0)=0$, and let $E$ be a closed  subset of $[0,2\pi]$  with $cap(E)>0$. Suppose that $f_\theta(x, h(x))$ is convergent for each $\theta\in E$. Then $f$ is convergent. \eT
\begin{pf} Let $g(x,y)=f((x+y)/2, -i(x-y)/2)$. Then $f(x,y)=g(x+iy, x-iy)$ and $f_\theta(x,y)=g(e^{i\theta}(x+iy), e^{-i\theta}(x-iy))$. Let $\phi_\theta(x)=f_\theta(x, h(x))=g(e^{i\theta}(x+ih(x)), e^{-i\theta}(x-ih(x)))$. Then $\phi_\theta(x)$  is convergent for $\theta\in E$. The $x$ terms of the two series $x\pm ih(x)$ cannot both be zero. Say, the $x$ term of $x+ih(x)$  is non-zero. So  $x+ih(x)$ has an inverse $\psi(x)$ which is a convergent power series such that $\psi(x)+ih(\psi(x))=x$. Set  $\psi(x)-ih(\psi(x))=\omega(x)$. Then $\phi_\theta(\psi(x))=g(e^{i\theta}x, e^{-i\theta}\omega(x))$ is convergent for $\theta\in E$. It follows that $g(sx,s^{-1}\omega(x))$ is convergent for each $s$ in the set $\{e^{i\theta}:\theta\in E\}$, which has positive capacity. By Theorem~\ref{MT1}, $g$ is convergent. Therefore $f$ is convergent.\end{pf}

\section{Analytic functions in ${\Bbb R}^2$}

Suppose that $f(x,y), \phi(x), q(x)$ are $C^\infty$ functions defined near the origin with $\phi(0)=0$.  Let $\hat f, \hat \phi, \hat q$ denote the Taylor series at 0 of those functions. Then $\hat f\in \CC[[x, y]]$, $\hat\phi,\hat q\in \CC[[x]]$. By the chain rule, $f(x,\phi(x))=q(x)$ implies $\hat f(x,\hat\phi(x))=\hat q(x)$. We consider here complex-valued analytic functions of real variables. If $I$ is an interval and if $\Gamma=\{(t,\gamma(t)): t\in I\}$ is a curve, the dilation by $s$ of $\Gamma$ is  $$\Gamma_s=\{(st, s\gamma(t))\}=\{(t, \gamma_{1/s}(t))\},\;\; \gamma_s(t)=s^{-1}\gamma(st).$$
\bT\label{FD} Let $f$ be a $C^\infty$ function defined in an open set $\Omega\subset\RR^2$ containing the origin, let $\Gamma=\{(t, \phi(t))\}$ be a non-linear analytic curve in $\RR^2$ passing through or ending at the origin, and let $E$ be a closed subset of ${\Bbb{R}}\backslash \{0\}$ of positive capacity. Suppose that for each $s\in E$, there is a real analytic function $F_s$ defined in a neighborhood $Q_s$ of $\Gamma_s\cap \Omega$ in $\RR^2$ such that $f$ and $F_s$ coincide on $\Gamma_s\cap \Omega$. Then there is a neighborhood $U$ of the origin, and a real analytic function $F$ defined on $U$ that coincides with $f$ on $U\cap \Lambda$, where $\Lambda:=\cup_{s\in E} \Gamma_s$.
\eT

\begin{pf} Without loss of generality we assume that $\phi(0)=0$. Since $f$ and $F_s$ coincide on $\Gamma_s$, we have
\begin{equation}\label{comp}f(x,\phi_{1/s}(x))=F_s(x,\phi_{1/s}(x)).\end{equation}
Let $g$, $h$ denote the Taylor series of $f$, $\phi$ respectively. Then (\ref{comp}) implies $g(x,h_{1/s}(x))=F_s(x,h_{1/s}(x))$. Hence $g(x,h_{1/s}(x))$ is convergent for $s\in E$. By Corollary~\ref{Dil}, $g$ is convergent. Thus $g$ represents a real analytic function $F$ in some neighborhood $U$ of the origin that satisfies $F(x,h_{1/s}(x))=F_s(x,h_{1/s}(x))$. It follows that the real analytic function $F$ coincides with $f$ on $U\cap\Lambda$. 
\end{pf}

Note that $f$ does not need to be analytic in a neighborhood of the origin.

If $\Gamma=\{(t,\phi(t):t\in I\}$ is a curve, its rotation by $\theta$ is 
$$\Gamma_\theta=\{(t\cos\theta+\phi(t)\sin\theta, -t\sin\theta+\phi(t)\cos\theta):t\in I\}.$$

\bT\label{FR} Let $f$ be a $C^\infty$ function defined in an open set $\Omega\subset\RR^2$ containing the origin, let $\Gamma=\{(t, \phi(t))\}$ be an analytic curve in $\RR^2$ passing through or ending at the origin, and let $E$ be a closed subset of $[0,2\pi]$ of positive capacity. Suppose that for each $\theta\in E$, there is a real analytic function $F_\theta$ defined in a neighborhood $Q_\theta$ of $\Gamma_\theta\cap \Omega$ in $\RR^2$ such that $f$ and $F_\theta$ coincide on $\Gamma_\theta\cap \Omega$. Then there is an analytic function $F$ defined in some neighborhood $U$ of the origin that coincides with $f$ on $U\cap \Lambda$, where $\Lambda:=\cup_{\theta\in E} \Gamma_\theta$. 
\eT
\begin{pf} The proof is similar to that of Theorem~\ref{FD}. Let $g_\theta(x,y):=g(x\cos\theta+y\sin\theta, -x\sin\theta+y\cos\theta)$. Then $g_\theta(x,h(x))$ is convergent for each $\theta\in E$. By Theorem~\ref{Rot}, $g$ is convergent.
\end{pf}

\bC\label{CFR} Let $f$ be a $C^\infty$ function defined in a neighborhood of \ $0$ in $\RR^2$, and let $\Gamma=\{(t, \phi(t))\}$ be an analytic curve passing through or ending at the origin in $\RR^2$. Suppose that for each $\theta\in[0,2\pi]$, the restriction of $f$ to $\Gamma_\theta$ extends to a real analytic function $F_\theta$ in a neighborhood $Q_\theta$ of  the origin. Then $f$ is analytic in a neighborhood of the origin.\eC

\noindent {\bf Remark.} We can see from the proofs that in Theorems 2.1, 2.2 and Corollary 2.3, the hypothesis on $f$ can be weakened to that $f$ has a Taylor series at the origin in the sense that there are numbers $a_{ij}$ such that for each positive integer $n$, 
$$f(x,y)-\sum_{i+j\le n} a_{ij}x^iy^j =o((x^2+y^2)^{n/2}).$$

\section{Examples}\label{exam} 
Here we present examples showing that the restrictions in our main theorems are necessary.

\vspace{5pt}
\noindent {\bf Example 3.1.} P. Lelong in the same paper \cite{Le} proved that if $E$ is a set with $cap(E)=0$ then one can find
a divergent power series $g(x,y)$ such that for all $s\in E$, $g(x,sx)$ is
convergent. For completeness we present here a construction of such an example.
Since $cap(E)=0$, there is a sequence of positive numbers $(\delta_n)$ with $\lim \delta_n=0$, and a sequence of polynomials $(P_n(x))$ with $\max_{x\in E} |P(x)|\le\delta_n^n$, where $P_n(x)=\sum_{j=0}^n b_{nj} x^{n-j}$ with $b_{n0}=1$. Let $a_{ij}=\delta_{i+j}^{-(i+j)}b_{i+j.i}$, and $g(x,y)=\sum a_{ij}x^iy^j$. Then $g(x,sx)=\sum \delta_n^{-n}P_n(s)x^n$. For $s\in E$, $|\delta_n^{-n}P_n(s)|\le1$, so $g(x,sx)$ is convergent. Note that $a_{0j}=\delta_j^{-j}$, which obviously implies that $g$ is divergent, since $\lim\delta_j=0$.\qed

\vspace{5pt} \noindent {\bf Example 3.2.} This example shows that the condition in Theorem~\ref{MT1} that $h(x)$ is not a monomial of the form $b_kx^k$ with $\sigma k-\tau=0$ cannot be dispensed with. Let $\sigma, k$ be positive integers, and $\phi\in\CC[[x]]$ a divergent series with $\phi(0)=0$. Let $g(x,y)=\phi(x^k)-\phi(y)$ and $h(x)=x^k$. Then $g$ is divergent; but $g(s^\sigma x, s^{\sigma k}h(x))=0$ for each $s\in\CC$.\qed

\vspace{5pt} \noindent {\bf Example 3.3.} This example shows that the hypothesis in Theorem~\ref{MT1} that $h(x)$ is convergent cannot be dispensed with when $\sigma\tau\le0$. (By Theorem~\ref{MT2} that hypothesis can be dispensed with when $\sigma\tau>0$.) The example also shows that Theorem~\ref{MT2} fails for $\sigma\tau\le0$.

Suppose that $\tau\le0$, $\sigma>0$. Let $u(x)=x+\cdots $ be a divergent series. Let $h(x), \phi(x)$ be the series satisfying $\phi(u(x))=x$ and $x^{|\tau|}h(x)^{\sigma}=u(x^{\sigma+|\tau|})$. Then $\phi, h$ are divergent. Let $f(x,y)=\phi(x^{|\tau|}y^\sigma)$. Then $f$ is divergent; but $f(s^\sigma x, s^{\tau}h(x))=x^{\sigma+|\tau|}$  for each $s\in\ {{\Bbb{C}}\backslash \{0\}}$. 

Now we consider the case where $\sigma=0$, $\tau=1$. Let $h(x)=x+\cdots $ be a divergent series, and let $\phi(x)$ be the series satisfying $h(x)\phi(x)=x^2$. Then $\phi$ is divergent. Let $f(x,y)=\phi(x)y$. Then $f$ is divergent; but $f(x, sh(x))=sx^2$ for each $s\in\CC$.\qed

We thank T.~S.~Neelon for informing us of the result of B.~Malgrange \cite{Ma} and pointing out several misprints. We are also grateful to the referee for his exhaustive and highly professional report. The changes implemented due to the suggestions in the report have significantly improved the paper.

\end{document}